\documentclass[11pt]{amsart}
\usepackage{amscd}

\theoremstyle{plain}
\newtheorem{theorem}{Theorem}%[section]
\newtheorem{corollary}[theorem]{Corollary}
\newtheorem{lemma}[theorem]{Lemma}
\newtheorem{proposition}[theorem]{Proposition}

\theoremstyle{definition}
\newtheorem{definition}[theorem]{Definition}
\newtheorem{remark}[theorem]{Remark}
\newtheorem*{questions}{Questions}
\newtheorem{example}[theorem]{Example}

\newcommand{\abs}[1]{\lvert#1\rvert}
\newcommand{\norm}[1]{\lVert#1\rVert}

\newcommand{\bignorm}[1]{\bigl\lVert#1\bigr\rVert}

\renewcommand{\le}{\leqslant}
\renewcommand{\ge}{\geqslant}\usepackage{amssymb}
\newcommand{\marg}[1]{}
\newcommand{\lab}[1]{\label{#1}\marg{#1}}
\renewcommand{\star}[1]{{}^*\!#1}

\newcommand{\term}[1]{{\textit{\textbf{#1}}}}

\def\mid{\::\:}
\def\re{r_{\rm ess}}
\def\se{\sigma_{\rm ess}}
\def\iE{i_{\scriptscriptstyle E}}
\def\iF{i_{\scriptscriptstyle F}}
\def\phiE{\varphi_{\scriptscriptstyle E}}
\def\phiF{\varphi_{\scriptscriptstyle F}}
\def\one{{\bf 1}}

\DeclareMathOperator{\ns}{ns}
\DeclareMathOperator{\fin}{fin}

\DeclareMathOperator{\dist}{dist}
\DeclareMathOperator{\co}{co}

\begin{document}

\input arrow.tex    % eplain package

\title[]{Measures of non-compactness\\
          of operators on Banach lattices}
\author[]{Vladimir G.~Troitsky}
\address{Department of Mathematical and Statistical Sciences, 632 CAB,
  University of Alberta, Edmonton, AB T6G\,2G1. Canada.}
\email{vtroitsky@math.ualberta.ca}

%\thanks{}
%\date{\today}

\begin{abstract}
  \cite{Andreu:91,Sadovsky:72} used representation spaces to study
  measures of non-compactness and spectral radii of operators on Banach
  lattices. In this paper, we develop representation spaces based on
  the nonstandard hull construction (which is equivalent to the
  ultrapower construction). As a particular application, we present a
  simple proof and some extensions of the main result
  of~\cite{dePagter:88} on the monotonicity of the measure of
  non-compactness and the spectral radius of AM-compact operators. We also
  use the representation spaces to characterize d-convergence and
  discuss the relationship between d-convergence and the measures of
  non-compactness.
\end{abstract}

\maketitle

\section{Introduction}

Recall that an operator $T$ between Banach lattices is said to be
\term{positive} if it maps positive vectors to positive vectors.  In
this case, we write $T\ge 0$. We write $S\le T$ if $T-S\ge 0$. We say
that $S$ is \term{dominated} by $T$ if $\abs{Sx}\le T\abs{x}$ for each
$x$.  An operator $T$ between Banach lattices is said to be
\term{order bounded} if it maps order intervals into order intervals.
It can be easily verified that if $T$ dominates $S$ then both $S$ and
$T$ are order bounded. An order bounded operator is \term{AM-compact}
if it maps order intervals (or almost order bounded sets) into
precompact sets. A set $A\subseteq E$ is \term{almost order bounded}
if for every $\varepsilon>0$ there exists $u\in E_+$ such that
$A\subseteq[-u,u]+\varepsilon B_E$, where $B_E$ stands for the unit
ball of $E$.  An operator $T$ between Banach lattices is said to be
\term{semicompact} if it maps bounded sets to almost order bounded
sets. We refer the reader
to~\cite{Aliprantis:85,Luxemburg:71,Meyer-Nieberg:91,Zaanen:83} for
a detailed study of Banach lattices and positive operators. All Banach
lattices in this paper are assumed to be complex unless specified
otherwise, all operators are assumed to be linear and bounded.

A lot of work has been done on the problem of the relationship between
compact operators and the order structure of a Banach lattice, see,
e.g.~\cite{Andreu:91,Schep:90,Wickstead:75,Wickstead:00}.  Still there
are many open questions. In particular, the problem can be considered
from the point of view of spectra of the operators.  It is well known
that if $T$ is a positive operator on a Banach lattice then the
spectral radius $r(T)$ belongs to the spectrum $\sigma(T)$. If $S$ is
dominated by $T$ then $\norm{S}\le\norm{T}$ and $r(S)\le r(T)$. The
central question of this paper is whether similar statements hold for
the essential spectrum, essential spectral radius, and the measure of
non-compactness.

Recall that the Calkin algebra of a Banach space $X$ is the quotient
of the algebra of $\mathcal L(X)$ over the closed algebraic ideal of
all compact operators. The \term{essential spectrum} $\se(T)$ and the
\term{essential spectral radius} $\re(T)$ of an operator $T$ on $X$
are defined as the spectrum and, respectively, the spectral radius of
the canonical image of $T$ in the Calkin algebra.

If $A$ is a bounded subset of a Banach space $X$ then the
\term{measure of noncompactness} $\chi(A)$ (it is sometimes referred
as the Hausdorff or ball measure of non-compactness) is defined via:
\begin{displaymath}
  \chi(A)=\inf\{\,\delta>0\mid A
  \mbox{ can be covered with a finite number
    of balls of radius }\delta\,\}.
\end{displaymath}
Clearly $\chi(A)=0$ if and only if $A$ is relatively compact. The
measure of noncompactness of an operator $T\colon X\to Y$ between
Banach spaces is defined via $\chi(T)=\chi(TB_X)$.  Then $\chi$ is a
seminorm on $\mathcal L(X,Y)$. It was shown in~\cite{Nussbaum:70}
that:
\begin{equation}
  \label{eq:nussbaum}\tag{*}
  \re(T)=\lim\limits_{n\to\infty}\sqrt[n]{\chi(T^n)}.
\end{equation}
for every $T\in\mathcal L(X)$. We refer the reader to~\cite{Banas:80}
for more details on measures of noncompactness.

\begin{questions}~
   \begin{enumerate}
   \item Does $\re(S)\in\se(S)$ for any positive operator $S$
     on a Banach lattice?
    \item Is $\re(S)\le\re(T)$ for any operators $S$ and
      $T$ provided $T$ dominates $S$?
      \lab{qi:main:re-re}
    \item Is $\chi(S)\le\chi(T)$ for any operators $T$ and
      $S$ provided $T$ dominates $S$?
      \lab{qi:main:chi-chi}
  \end{enumerate}
\end{questions}

These questions were first addressed in~\cite{dePagter:88}, and the
following results were obtained (see also
\cite[Section 4.3]{Meyer-Nieberg:91}):

\begin{theorem}[{\cite{dePagter:88}}]\lab{t:dPS-main}~
   \begin{enumerate}
    \item If $S$ is a positive AM-compact operator on a Banach lattice 
      then $\re(S)\in\se(S)$. \lab{ti:dPS-main:re-se}
    \item If $S$ and $T$ are two operators on a Banach lattice such
      that $0\le S\le T$ and $S$ is AM-compact, then
      $\re(S)\le\re(T)$.\lab{ti:dPS-main:re-re} 
    \item If $S,T\colon E\to F$ are two operators between Banach
      lattices such that $0\le S\le T$, $S$ is AM-compact, and both
      $E'$ and $F$ have order continuous norms, then
      $\chi(S)\le\chi(T)$. \lab{ti:dPS-main:chi-chi}
  \end{enumerate}
\end{theorem}

In the same paper, an example is given of two
non-AM-compact operators $0\le S\le T$ with $\re(S)>\re(T)$ and
$\re(S)\notin\se(S)$.

It is easy to see that if an operator $S\colon E\to E$ is dominated by
a compact operator $T$ then $\re(S)=0$. Indeed, $0\le S+T\le 2T$, so
that $(S+T)^3$ is compact by the Cube Theorem \cite[Theorem
16.14]{Aliprantis:85}, and it follows that $\re(S)=0$. But, as far as
we know, it is still not known whether every operator dominated by an
essentially quasinilpotent operator is itself essentially
quasinilpotent.

An important technical tool used in~\cite{dePagter:88} is
a \term{measure of non-semicompactness}, introduced analogously to the
Hausdorff measure of noncompactness: if $A$ is a norm bounded set in a
Banach lattice, then
\begin{displaymath}
  \rho(A)=\inf\{\,\delta>0 \mid \,
    A\subseteq [-u,u]+\delta B_E\mbox{ for some }u\in E_+\,\},
\end{displaymath}
and $\rho(T)=\rho(TB_E)$ whenever $T$ is an operator between Banach
lattices.  Clearly, $\rho(A)=0$ iff $A$ is almost order bounded, and
$\rho(T)=0$ iff $T$ is semi-compact. Furthermore, $\rho(A)\le\chi(A)$ for
each bounded set $A$ and $\rho(T)\le\chi(T)$ for each bounded operator
~$T$.  It was proved in~\cite[Theorem~2.5]{dePagter:88} that
$\rho(T)=\chi(T)$ for every order bounded AM-compact operator $T\colon
E\to F$ when $E'$ and $F$ have order continuous norms.

In this paper, we develop a certain representation space technique and,
using this technique, we present a simple proof of
Theorem~\ref{t:dPS-main} as well as some improvements of it.  Our
technique is based on the nonstandard hull construction of Nonstandard
Analysis (which is equivalent to the ultrapower construction). Recall
the construction briefly.  If~$X$ is a (standard) Banach space, denote
by $\star{X}$ the \term{nonstandard extension} of~$X$. The symbols
$x$, $y$, $z$, etc., will usually stand for elements of~$\star{X}$.
If~$T$ is a standard operator on~$X$, we use the same symbol~$T$
instead of~$\star{T}$ to denote the extension of $T$ to~$\star{X}$.
The symbol $\fin(\star{X})$ stands for the subspace of all elements of
$\star{X}$ of finite norm, while the \term{monad of zero}, $\mu(0)$,
consists of those elements of $\star{X}$ whose norm is infinitesimal.
The Banach space $\widehat{X}=\fin(\star{X})/\mu(0)$ is called the
\term{nonstandard hull} of $X$. If $x\in\fin(\star{X})$ then $\hat{x}$
will stand for the corresponding element in $\widehat{X}$. Every
bounded operator between Banach spaces $T\colon X\to Y$ induces a
bounded operators $\widehat{T}\colon\widehat{X}\to\widehat{Y}$ via
$\widehat{T}\hat{x}=\widehat{Tx}$.  Clearly, $X$ is isometrically
isomorphic to $\ns(\star{X})/\mu(0)$, where $\ns(\star{X})$ stands for
the set of all near-standard elements of $\star{X}$.  Thus, one can
view $X$ as a closed subspace of~$\widehat{X}$.  For $A\subseteq X$ we
write $\widehat{A}=\{\hat{x}\mid x\in\star{A}\}$. It is known that a
standard set $A\subseteq X$ is relatively compact iff
$\star{A}\subseteq\ns(\star{X})$ iff $\widehat{A}\subseteq X$.  If $E$
is a Banach lattice then $\widehat{E}$ also is a Banach lattice.
Further details on nonstandard analysis and nonstandard hulls can be
found in~\cite{Cozart:74,Henson:83,Wolff:97}.
A reader familiar with the technique of ultraproducts can view
$\widehat{X}$ as an ultrapower of~$X$. Clearly, all the proofs in this
paper can be redone in terms of ultrapowers, but we believe that the
language of Nonstandard Analysis is more appropriate for this problem.

Define $\widetilde{X}=\widehat{X}/X$. We will see in
Section~\ref{sec:rep-spaces} that this space is a representation space
for $\chi$ and $\se$, that is, $\chi(T)=\norm{\widetilde{T}}$ and
$\se(T)=\sigma(\widetilde{T})$. For a Banach lattice $E$, let $I(E)$
be the order ideal generated by $E$ in $\widehat{E}$, and let
$\iE=\overline{I(E)}$. Now we define $\check{E}=\widehat{E}/\iE$.  It
will be shown in Section~\ref{sec:rep-spaces} that $\check{E}$ is a
representation space for $\rho(T)$, that is,
$\rho(T)=\norm{\check{T}}$.  We will also show that $T\colon E\to F$
is AM-compact if and only if $\widehat{T}$ maps $\iE$ into $F$.  In
Section~\ref{sec:appl} we use these representation spaces to prove
Theorem~\ref{t:dPS-main} and similar results. It should be mentioned
that various representation spaces have been used to study the essential
spectrum of an operator, see e.g.,~\cite{Andreu:91,Sadovsky:72}.

Finally, in Section~\ref{sec:d-top} we discuss the d-topology on a
Banach lattice. We say that a net $(x_\alpha)$ \term{d-converges} to
$x$ in $E$ if $\abs{x_\alpha-x}\wedge y$ converges to zero in norm for
every $y\in E_+$. We investigate the relation between the d-topology,
the space $\widehat{E}$, and the questions stated in the beginning of
the paper. We also discuss examples of d-topologies.

\section{Representation Spaces} \label{sec:rep-spaces}

Given a Banach space $X$, we define
$\widetilde{X}=\widehat{X}/X$. If $x\in \star{X}$, then $\tilde{x}$
will stand for the corresponding element in~$\widetilde{X}$. Every
operator $T\colon X\to Y$ induces an operators
$\widetilde{T}\colon\widetilde{X}\to\widetilde{Y}$ given by
$\widetilde{T}\tilde{x}=\widetilde{Tx}$. An
operator $T\colon X\to Y$ is compact iff the range of $\widehat{T}$ is
contained in~$Y$ iff $\widetilde{T}=0$.

\begin{lemma} \lab{l:chi-of-set}
  If $A$ is a bounded subset of $X$ then
  $\chi(A)=\max\limits_{y\in\star{A}}\norm{\tilde{y}}$.
\end{lemma}

\begin{proof}
  Fix $y\in\star{A}$. For every standard $\varepsilon>0$ there is a
  finite set $F\subset X$ such that $A$ (and hence $\star{A}$) is within
  $\chi(A)+\varepsilon$ of $F$. Then
  $\norm{\tilde{y}}\le\dist(\hat{y},F)\le\chi(A)+\varepsilon$, so that
  $\norm{\tilde{y}}\le\chi(A)$.  Conversely, for every finite family
  of balls of radius less than $\chi(A)$ there exists a point in $A$
  which is not covered by the balls. By the idealization (saturation)
  principle, there exists $y\in\star{A}$ which does not belong to any
  standard ball of radius less than $\chi(A)$. Therefore
  $\norm{\tilde{y}}\ge\chi(A)$.
\end{proof}

From this point of view, $\chi(A)$ measures how far the set $\star{A}$
is from $\ns(\star{X})$.  Then $\chi(T)$ measures how much closer sets
become to $\ns(\star{X})$ after we apply $T$. Actually, the following
lemma describes the relation. Similar results were proved
in~\cite{Wolff:97,Sadovsky:72}.

\begin{corollary} \lab{c:chi-of-operator}
  If $T$ is a bounded operator on $X$ then $\chi(T)=\norm{\widetilde{T}}$.
\end{corollary}

Corollary~\ref{c:chi-of-operator} and
formula~(\ref{eq:nussbaum}) imply that $r_{\rm{ess}}(T)=r(\widetilde{T})$
for every bounded operator on a Banach space. Moreover, it was shown
in~\cite{Buoni:77} and in~\cite[Theorem 3.11]{Wolff:97} that
$\widetilde{T}$ is invertible iff $T$ is Fredholm, so that
$\sigma(\widetilde{T})=\sigma_{\rm ess}(T)$.

\medskip

Now we turn to Banach lattices.  It is well known that if $E$ is a
Banach lattice then $\widehat{E}$ is also a Banach lattice.

\begin{remark}\lab{r:intervals}
  The following simple observation turns out to be quite handy in the
  context of our problem. Suppose that $I$ is an (order) ideal in a
  vector lattice $E$ and consider the quotient vector lattice $E/I$.
  It is known that the canonical epimorphism from $E$ onto $E/I$ is a
  lattice homomorphism, hence it maps order intervals onto order
  intervals. It follows that if $a,b,x\in E$ are such that $a\le b$
  and $[a]\le[x]\le[b]$, where $[a]$, $[x]$, and $[b]$ are the
  equivalence classes of $a$, $x$, and $b$ respectively in the
  quotient vector lattice $E/I$, then $a\le x'\le b$ for some
  $x'\in[x]$.  In particular, if $E$ is a Banach lattice and
  $\hat{a}\le\hat{x}\le\hat{b}$ for some $a,x,b\in\star{E}$ such that
  $a\le b$, then $a\le x'\le b$ for some $x'\in\star{E}$ such that
  $x'\approx x$.
\end{remark}

The following important characterization was obtained
in~\cite{Cozart:74}.

\begin{theorem} \lab{t:cm-hull}
  The following statements are equivalent:
  \begin{enumerate}
    \item $\widehat{E}$ is Dedekind complete;
    \item $\widehat{E}$ has the projection property;
    \item $\widehat{E}$ has order continuous norm;
    \item $c_0$ is not lattice finitely representable%
      \footnote{Recall that a Banach lattice $F$ is lattice finitely
        representable in a Banach lattice $E$ if for each finite
        dimensional vector sublattice $H$ in $F$ and for each
        $\delta>0$ there exists a lattice embedding $T\colon H\to E$
        such that $\norm{T},\norm{T^{-1}}\le 1+\delta$.} 
      in $E$.
  \end{enumerate}
\end{theorem}

Notice that $\widetilde{E}$ need not be a Banach lattice because $E$
might not be an (order) ideal in $\widehat{E}$. In fact, this happens
only when $E$ is atomic with order continuous norm.
It was first noticed in~\cite{Cozart:74} that $E$ is an ideal in
$\widehat{E}$ if and only if the order intervals in $E$ are compact.
Indeed, $[-u,u]$ is a compact set in $E$ for each $u\in E_+$ if and only
if $[-u,u]_{\star{E}}\subset\ns(\star{E})$. In view of
Remark~\ref{r:intervals} this is equivalent to
$[-\hat{u},\hat{u}]_{\widehat{E}}\subset E$ for each $u\in E_+$.
It was shown in~\cite{Wickstead:75,Walsh:68} that the intervals in a
Banach lattice $E$ are compact if and only if $E$ is atomic with order
continuous norm. Thus, the following result holds.

\begin{proposition}  \lab{p:discrete-lattice}
  The space $\widetilde{E}$ is a Banach lattice if and only if $E$ is
  atomic with order continuous norm. 
\end{proposition}

Denote by $I(E)$ the ideal generated by $E$ in $\widehat{E}$. This
ideal was extensively studied in~\cite{Cozart:74}.

\begin{theorem}[{\cite{Cozart:74}}] \lab{t:cm-ideal}
  The following statements are equivalent:
  \begin{enumerate}
    \item $I(E)$ is Dedekind complete;
    \item $I(E)$ has the projection property;
    \item $I(E)$ has order continuous norm;
    \item $E$ has order continuous norm;
  \end{enumerate}
\end{theorem}

It also follows from Proposition~\ref{p:discrete-lattice} that $I(E)=E$ if
and only if $E$ is atomic with order continuous norm.

\begin{example}
  {\it $I(E)$ need not be norm closed.} Let
  $E=L_1[0,1]$, fix an infinite positive integer $N$, and consider the
  partition of $\star{[0,1]}$ into $2N$ equal subintervals. Let $A$ be
  the union of all odd-numbered intervals, i.e.,
  $A=\bigcup_{i=1}^N[\frac{i}{N}-\frac{1}{2N},\frac{i}{N}]$, and
  $f=\chi_A$, the characteristic function of $A$. We claim that zero
  is the greatest standard function in $\widehat{E}$ dominated by
  $\hat{f}$.  Indeed, suppose that $\hat{g}\le\hat{f}$ for some $g\in
  L_1[0,1]$ such that $g$ is positive on a set of positive measure.
  Then there exists $\varepsilon>0$ such that $m(C)>0$, where
  $C=\{g>\varepsilon\}$ and $m$ stands for the Lebesgue measure.
  Clearly $(\varepsilon\chi_{C}-f)^+\le(g-f)^+\approx 0$. Let
  $\delta=\frac{\varepsilon\cdot m(C)}{4+\varepsilon}$. Then one can
  find a set $D\subseteq[0,1]$ such that $D$ is a finite union of
  intervals and $m(C\bigtriangleup D)<\delta$. It follows that
  $\bignorm{(\varepsilon\chi_D-f)^+}\le
  \bignorm{(\varepsilon\chi_{C}-f)^+}+\delta\lesssim\delta$.  On the
  other hand, since $m(A\cap I)=\frac{1}{2}m(I)$ for every standard
  interval $I\subseteq[0,1]$, we have
  $\bignorm{(\varepsilon\chi_D-f)^+}=\frac{\varepsilon}{2}m(D)\ge
  \frac{\varepsilon}{2}\bigl(m(C)-\delta\bigr)= 2\delta$, a
  contradiction. It can be shown in a similar fashion that
  $\chi_{[0,1]}$ is the smallest standard function that dominates $f$.
  
  Now let $E=L_1(\mathbb R)$. Again, fix an infinite positive integer
  $N$, and let $A_1$ be the set $A$ from the previous paragraph. Cut
  the interval $\star{[1,2]}$ into $4N$ equal subintervals and let
  $A_2$ be the union of every fourth interval, i.e.,
  $A_2=1+\bigcup_{i=1}^N[\frac{i}{N}-\frac{1}{4N},\frac{i}{N}]$.
  Similarly, for every $n\in\star{\mathbb N}$ let
  $A_n=n-1+\bigcup_{i=1}^N[\frac{i}{N}-\frac{1}{2^nN},\frac{i}{N}]$.
  Next, let $B_n=\bigcup_{k=1}^nA_k$ for each $n\in\star{\mathbb N}$
  and $B=\bigcup_{k\in\star{\mathbb N}}A_k$. Then
  $m(B)=\frac{1}{2}+\frac{1}{4}+\frac{1}{8}+\dots<\infty$.  For each
  $n\in\star{\mathbb N}$ let $h_n$ be the characteristic function of
  $B_n$, and let $h$ be the characteristic function of $B$.  Notice
  that $\norm{h-h_n}=m(B\setminus
  B_n)=\frac{1}{2^{n+1}}+\frac{1}{2^{n+2}}+\dots$, so that
  $\norm{\hat{h}-\hat{h}_n}\to 0$ as $n\to 0$ in~$\mathbb N$.  On
  the other hand, $h_n\le\chi_{[0,n]}$ for every $n$, so that
  $\hat{h}_n\in I(E)$ for every standard~$n$. But it follows from the
  previous paragraph, that $h$ is not dominated by a standard function
  in $L_1(\mathbb R)$, because this function would have to be greater
  or equal than 1 a.e. on $\mathbb R$. Therefore, $\hat{h}$ is in the
  closure of $I(E)$ but not in $I(E)$, so that $I(E)$ is not closed.
\end{example}

Denote by $\iE$ the closure of $I(E)$ in $\widehat{E}$. 

\begin{remark} \lab{r:iE-projband}
  Since every closed ideal in a Banach lattice with order continuous
  norm is a band, it follows from Theorem~\ref{t:cm-hull}
  that if $c_0$ is not lattice finitely representable in $E$ then
  $\iE$ is a projection band.
\end{remark}

Define a representation space $\check{E}$ for $E$ via
$\check{E}=\widehat{E}/\iE$. Apparently $\check{E}$ is a Banach
lattice. It follows from Lemma~\ref{p:discrete-lattice} that
$\check{E}=\widetilde{E}$ if and only if $E$ is atomic with order
continuous norm.  On the other hand, if $E$ has a strong order unit,
then $I(E)=\widehat{E}$ and $\check{E}$ is trivial. It follows from
$E\subseteq\iE\subseteq\widehat{E}$ that
$\norm{\check{x}}\le\norm{\tilde{x}}\le\norm{\hat{x}}$ for each
$x\in\fin(\star{E})$.

Let $T\colon E\to F$ be an order bounded operator between Banach
lattices. Suppose that $\hat{x}\in I(E)$. Then $\abs{\hat{x}}\le\hat{u}$
for some $u\in E_+$. By Remark~\ref{r:intervals} there exists
$y\in\star{E}$ such that $x\approx y\in\star{[-u,u]}$. Then $Tx\approx Ty\in
T\star{[-u,u]}\subseteq\star{[-v,v]}$ for some $v\in F_+$, so that
$\widehat{T}\hat{x}\in I(F)$. Thus, $\widehat{T}$ maps $\iE$ into
$\iF$ and, therefore, $\check{x}\mapsto(Tx)\check{}$ induces a bounded
operator $\check{T}$ from $\check{E}$ to $\check{F}$.

We claim that $\check{T}$ is order bounded.  Indeed, suppose that
$\check{u}>0$ and $\check{x}\in[-\check{u},\check{u}]$ for some
$x,u\in\fin(\star{E})$. By Remark~\ref{r:intervals} we can assume that
$u\in\star{E}_+$ and $x\in[-u,u]$. Then $Tx\in T[-u,u]\subseteq[-v,v]$
for some $v\in\star{E}_+$. Therefore,
$\check{T}[-\check{u},\check{u}]\subseteq[-\check{v},\check{v}]$, hence
$\check{T}$ is order bounded.  Notice that $v$ can be chosen in
$\fin(\star{E})$ because of the following fact, which is due to
A.~Wickstead: {\it if $T\colon E\to F$ is order bounded then}
\begin{displaymath}
    \sup\limits_{\,\norm{z}\le 1,z\ge 0}\inf\bigl\{\norm{y}\mid y
    \in E_+\mbox{ and }
    T[-z,z]\subseteq[-y,y]\,\bigr\}\quad<\quad\infty. 
\end{displaymath}
Indeed, otherwise for each $n>0$ one could find a positive $z_n$ in
$B_E$ with $$\inf\bigl\{\norm{y}\mid y\in E_+\mbox{ and }T[-z_n,z_n]
\subseteq[-y,y]\bigr\}>n^3.$$ Let $z=\sum_{n=1}^\infty\frac{z_n}{n^2}$,
it is easy to see that $T[-z,z]$ is not contained in any order
interval in~$F$.

Clearly, if $T$ is positive then $\check{T}$ is positive and
$\check{T}=0$ if and only if $T$ is semi-compact.

The following two results are analogous to Lemma~\ref{l:chi-of-set}
and Corollary~\ref{c:chi-of-operator}.

\begin{lemma}  \lab{l:rho-of-set}
  If $A$ is a bounded subset of $E$, then
  $\rho(A)=\max\limits_{y\in\star{A}}\norm{\check{y}}$.  
\end{lemma}

\begin{proof}
  Fix a standard $\gamma>\rho(A)$. Then $A\subseteq[-u,u]+\gamma B_E$
  for some $u\in E_+$. For each $y\in\star{A}$ we have $y=v+h$
  such that $v\in\star{[-u,u]}$ and $\norm{h}\le\gamma$. It follows
  from $\hat{v}\in\iE$ that $\norm{\check{y}}=\norm{\hat{h}}\le\gamma$. Thus
  $\norm{\check{y}}\le\rho(A)$.
  
  On the other hand, for every standard positive $\gamma<\rho(A)$ and
  for every $u\in E_+$ there is a point $y$ in $A$ which doesn't
  belong to $[-u,u]+\gamma B_E$. By saturation there exists
  $y\in\star{A}$ such that $y\notin\star{[-u,u]}+\gamma\star{B_E}$ for every
  standard positive $\gamma<\rho(A)$ and for every $u\in E_+$. Then
  $\norm{\check{y}}\ge\rho(A)$.
\end{proof}

It follows, in particular (c.f.~\cite[Corollary 1.4]{Schep:90}), that
a bounded subset $A\subset E$ is almost order bounded if and only if
$\widehat{A}\subseteq\iE$.

\begin{corollary} \lab{c:rho-of-operator}
  If $T\colon E\to F$ is operator between Banach lattices then
  $\rho(T)=\norm{\check{T}}$.
\end{corollary}

\section{Applications} \label{sec:appl}

The following theorem follows immediately from
Proposition~\ref{p:discrete-lattice} and
Corollary~\ref{c:chi-of-operator}.

\begin{theorem}
  If $E$ and $F$ are atomic Banach lattices with
  order continuous norms then:
   \begin{enumerate}
    \item if $T$ is a positive operator on $E$ then $\re(T)\in\se(T)$;
    \item if $S,T\colon E\to E$ and $T$ dominates $S$ then
      $\re(S)\le\re(T)$;
    \item \lab{ti:discr:chi-chi} if $S,T\colon E\to F$ and $T$
      dominates $S$ then $\chi(S)\le\chi(T)$.
  \end{enumerate}
\end{theorem}

Since every operator on an atomic Banach lattice is AM-compact, this
result can be viewed as a special case of Theorem~\ref{t:dPS-main}
except that we do not require $E'$ to have order continuous norm
in~(\ref{ti:discr:chi-chi}).

\medskip

Next, we are going to characterize AM-compact operators.
Denote by $\phiE$ the canonical epimorphism from
$\widetilde{E}$ to $\check{E}$ given by $\phiE(\tilde{x})=\check{x}$.
By the definition of $\check{T}$ we have
$\check{T}\phiE=\phiF\widetilde{T}$.

\begin{theorem}  \lab{t:am}
  Let $T\colon E\to F$ be an order bounded operator between Banach
  lattices. The following statements are equivalent:
  \begin{enumerate}
    \item \lab{ti:am}    $T$ is AM-compact;
    \item \lab{ti:hat}   $\widehat{T}$ maps $\iE$ into $F$;
    \item \lab{ti:diag}  There exists a map
      $\bar{T}\colon\check{E}\to\widetilde{F}$ such that
      $\widetilde{T}=\bar{T}\phiE$,
      i.e., $\widetilde{T}\tilde{x}=\bar{T}\check{x}$:
      \begin{displaymath}
        \commdiag{\widetilde{E}&\mapright^{\widetilde{T}}&\widetilde{F}\cr
                  \mapdown\lft{\varphi_E}&\arrow(3,2)\rt{\bar{T}}\cr
                  \check{E}}
      \end{displaymath}
  \end{enumerate}
\end{theorem}

\begin{proof}
  If $u\in E_+$ then $T[-u,u]$ is compact if and only if
  $\widehat{T[-u,u]}\subseteq F$. In view of Remark~\ref{r:intervals}
  we have $\widehat{T[-u,u]}=\widehat{T}[-\hat{u},\hat{u}]$, so that
  (\ref{ti:am})$\Leftrightarrow$(\ref{ti:hat}). To show
  (\ref{ti:hat})$\Leftrightarrow$(\ref{ti:diag}) notice that
  $\ker\phiE=\iE/E$. If $\widehat{T}$ maps $\iE$ into $F$ then
  $\iE/E\subseteq\ker\widetilde{T}$, so that
  $\bar{T}\bigl(\phiE(\tilde{x})\bigr)=\widetilde{T}\tilde{x}$
  defines an operator from $\check{E}$ to $\widetilde{F}$.
  Conversely, if such a $\bar{T}$ exists, then for every
  $\hat{x}\in\iE$ we have $\phiE(\tilde{x})=0$ so that
  $\widetilde{T}\tilde{x}=\bar{T}\phiE(\tilde{x})=0$, hence
  $\widehat{T}\hat{x}\in F$.
\end{proof}

\begin{remark} \lab{r:product}
  Notice that if $T$ is AM-compact then
  $\norm{\widetilde{T}}=\norm{\bar{T}}$ because $\phiE$ maps the unit
  ball of $\widetilde{E}$ onto the unit ball of $\check{E}$.

  If $E$, $F$, and $G$ are Banach lattices and $S\colon E\to
  F$ and $T\colon F\to G$ are order bounded operators such that $T$ is 
  AM-compact, then it follows from the diagram
  \begin{displaymath}
    \commdiag{
      \widetilde{E}&\mapright^{\widetilde{S}}&\widetilde{F}&
                        \mapright^{\widetilde{T}}&\widetilde{G}\cr
      \mapdown\lft{\phiE}&&\mapdown\lft{\phiF}&\arrow(3,2)\lft{\bar{T}}&\cr
      \check{E}&\mapright^{\check{S}}&\check{F}&&}
  \end{displaymath}
  that $\norm{\widetilde{T}\widetilde{S}}\le
  \norm{\phiE}\norm{\check{S}}\norm{\bar{T}}=
  \norm{\check{S}}\norm{\widetilde{T}}$ (c.f.
  \cite[Lemma~3.1]{dePagter:88}.
\end{remark}

Now we are ready to present a simple proof of
Theorem~\ref{t:dPS-main}. We replace the condition $0\le S\le T$ with
the slightly weaker condition of $S$ being dominated by $T$.

\begin{proof}[\bf Proof of Theorem~\ref{t:dPS-main}]
  If $T$ is positive and AM-compact, then by Remark~\ref{r:product} we
  have $\norm{\widetilde{T}^n}\le
  \norm{\check{T}^{n-1}}\norm{\widetilde{T}}$, so that
  \begin{displaymath}
    r(\widetilde{T})=
    \lim\limits_{n\to\infty}\sqrt[n]{\norm{\widetilde{T}^n}}\le
    \lim\limits_{n\to\infty}\sqrt[n]{\norm{\check{T}^{n-1}}\norm{\widetilde{T}}}=
    \lim\limits_{n\to\infty}\sqrt[n]{\norm{\check{T}^{n-1}}}=r(\check{T}).
  \end{displaymath}
  On the other hand, we always have $r(\check{T})\le
  r(\widetilde{T})$, so that $\re(T)=r(\widetilde{T})=r(\check{T})$.
  Since $\check{T}$ is a positive operator on the Banach lattice
  $\check{E}$ we have $r(\check{T})\in\sigma(\check{T})$. Finally,
  $\check{T}$ is a quotient of $\widetilde{T}$ so that
  $\sigma(\check{T})\subseteq\sigma(\widetilde{T})=\sigma_{\rm
    ess}(T)$. Thus, $\re(T)\in\sigma_{\rm ess}(T)$.

  Next, if $T$ dominates $S$ and $S$ is AM-compact then by
  Remark~\ref{r:product} we have
  $\norm{\widetilde{S}^n}\le\norm{\check{S}^{n-1}}\norm{\widetilde{S}}
  \le\norm{\check{T}^{n-1}}\norm{\widetilde{S}}$. This yields
  \begin{displaymath}
    \re(S)=r(\widetilde{S})=
    \lim\limits_{n\to\infty}\sqrt[n]{\norm{\widetilde{S}^n}}\le
    \lim\limits_{n\to\infty}\sqrt[n]{\norm{\check{T}^{n-1}}}=
    r(\check{T})\le r(\widetilde{T})=\re(T).
  \end{displaymath}
  
  If in addition $E'$ and $F$ have order continuous norms, then
  $\chi(S)=\rho(S)=\norm{\check{S}}$ and
  $\norm{\check{T}}=\rho(T)\le\chi(T)$, but since $\check{E}$ is a
  Banach lattice and $\check{T}$ dominates $\check{S}$ we conclude
  that $\norm{\check{S}}\le\norm{\check{T}}$ so that
  $\chi(S)\le\chi(T)$.
\end{proof}

The following theorem is an analog of~\cite[Theorem~1.5]{Andreu:91}
for $\widetilde{T}$ and $\check{T}$.

\begin{theorem} \lab{t:am-spectra}
  Let $E$ be a Banach lattice and $T\colon E\to E$ an AM-compact
  operator. Then
  \begin{enumerate}
    \item \lab{ti:pspectra} $\sigma_p(\check{T})\setminus\{0\}=
            \sigma_p(\widetilde{T})\setminus\{0\}$;
    \item \lab{ti:onto} If $\lambda\neq 0$ then $\lambda
      I-\check{T}$ is onto if and only if $\lambda I-\widetilde{T}$ is
      onto;
    \item \lab{ti:spectra} $\sigma(\check{T})\setminus\{0\}=
            \sigma(\widetilde{T})\setminus\{0\}=
            \sigma_{\rm ess}(T)\setminus\{0\}$.
  \end{enumerate}
\end{theorem}

\begin{proof}
  (\ref{ti:pspectra}) Suppose that $\lambda\neq 0$. If
  $\check{T}\check{x}=\lambda\check{x}$, $\check{x}\neq 0$ then
  $(\lambda I-\widehat{T})\hat{x}\in\iE$ so that $\widehat{T}(\lambda
  I-\widehat{T})\hat{x}\in E$ hence $(\lambda
  I-\widetilde{T})(\widetilde{T}\tilde{x})=0$. Notice that
  $\widetilde{T}\tilde{x}\neq 0$ because
  $\check{T}\check{x}=\lambda\check{x}\neq 0$. This yields
  $\lambda\in\sigma_p(\widetilde{T})$. Conversely, if
  $\widetilde{T}\widetilde{x}=\lambda\widetilde{x}$,
  $\widetilde{x}\neq 0$ then $\check{T}\check{x}=\lambda\check{x}$.
  Notice that $\check{x}\neq 0$ because otherwise we would have
  $\hat{x}\in\iE$ which would imply $\widehat{T}\hat{x}\in E$ and
  $\widetilde{T}\tilde{x}=0$.
  
  (\ref{ti:onto}) If $\lambda I-\widetilde{T}$ is onto, then $\lambda
  I-\check{T}$ is also onto as a quotient of $\lambda
  I-\widetilde{T}$.  To show the converse, take $\hat{y}\in\widehat{E}$,
  then there exists $\hat{x}\in\widehat{E}$, such that $(\lambda
  I-\check{T})\check{x}=\check{y}$. Then $(\lambda
  I-\widehat{T})\hat{x}-\hat{y}\in\iE$ so that
  $\widehat{T}\bigl((\lambda I-\widehat{T})\hat{x}-\hat{y}\bigr)\in
  E$. This yields $(\lambda
  I-\widetilde{T})(\widetilde{T}\tilde{x}+\tilde{y})-\lambda\tilde{y}=0$.
  Hence $(\lambda I-\widetilde{T})(\frac{1}{\lambda})
  (\widetilde{T}\tilde{x}+\tilde{y})=\tilde{y}$ and, therefore,
  $\lambda I-\widetilde{T}$ is onto.
  
  Finally, (\ref{ti:spectra}) follows immediately from
  (\ref{ti:pspectra}) and (\ref{ti:onto}).
\end{proof}

\section{d-topologies} \label{sec:d-top}

\begin{definition}
  We say that a net $(x_\alpha)$ in a Banach lattice $E$
  \term{d-converges} to $x\in E$ if $\abs{x_\alpha-x}\wedge y$
  converges to zero in norm for every $y\in E_+$. The topology
  generated by this convergence will be referred to as the
  \term{d-topology} of $E$.
\end{definition}

Let $\mu_d(0)$ denotes the monad of zero for the d-topology, while $E^d$
stands for the disjoint complement of $E$ in $\widehat{E}$.

\begin{lemma} \lab{l:mud}
  For a point $x\in\star{E}$ the following are equivalent:
  \begin{enumerate}
    \item \lab{li:mu}  $x\in\mu_d(0)$;
    \item \lab{li:hatd} $\hat{x}\in E^d$;
    \item \lab{li:nd}  $x$ is \term{nearly disjoint} with every $y\in E_+$,
      i.e., $\abs{x}\wedge y\approx 0$.
  \end{enumerate}
\end{lemma}

\begin{proof}
  To see that (\ref{li:hatd})$\Leftrightarrow$(\ref{li:nd}) observe
  that $\hat{x}\in E^d$ if and only if for every $y\in E_+$ we have
  $\abs{\hat{x}}\wedge\hat{y}=0$ or, equivalently, $\abs{x}\wedge
  y\approx 0$. On the other hand, $(x_\alpha)$ d-converges to zero in
  $E$ if and only if $\abs{x_\alpha}\wedge y$ converges to zero in
  norm for every $y\in E_+$, which is equivalent to
  $\abs{x_\alpha}\wedge y\approx 0$ for every infinite $\alpha$, so
  that (\ref{li:mu})$\Leftrightarrow$(\ref{li:nd}).
\end{proof}

Notice that the embedding of $E^d$ into $\check{E}$ given by
$\hat{x}\mapsto\check{x}$ is an isometry (we will see in
Example~\ref{e:c-0} that it need not be onto). Indeed, if $\hat{x}\in
E^d$ then clearly $\norm{\check{x}}\le\norm{\hat{x}}$.  On the other
hand, $\hat{x}\perp\hat{y}$ for each $\hat{y}\in\iE$, so that
$\norm{\hat{x}+\hat{y}}\ge\norm{\hat{x}}$, whence
$\norm{\check{x}}=\norm{\tilde{x}}$.  It is easy to see that
$\norm{\check{x}}=\norm{\tilde{x}}$ for every $\hat{x}\in E\oplus
E^d$. It follows then from Lemmas~\ref{l:chi-of-set}
and~\ref{l:rho-of-set} that if $A$ is a bounded set in $A$ then
$\chi(A)=\rho(A)$ whenever $\widehat{A}\subseteq E\oplus E^d$. By
Lemma~\ref{l:mud} $\hat{x}\in E\oplus E^d$ means that $x$ is
nearstandard relative to the d-topology, so that {\it $\widehat{A}\subseteq
E\oplus E^d$ if and only if $A$ is relatively d-compact\/} (i.e.,
relatively compact with respect to the d-topology).  Thus, we arrive
to the following result.

\begin{proposition} \lab{p:dcomp-set-chi-rho}
  If $A\subseteq E$ is d-compact, then $\chi(A)=\rho(A)$.
\end{proposition}

We say that an operator between Banach lattices
is \term{d-compact} if it maps bounded sets into relatively d-compact
sets. Observe that if $T$ is order bounded and d-compact then it is
AM-compact. Indeed, if $y\in E_+$ then $T[-y,y]$ is relatively
d-compact and order bounded, but the d-topology agrees with the norm
topology on order bounded sets, so that $T[-y,y]$ is relatively
compact. It follows, in particular, that Theorem~\ref{t:dPS-main}
applies to d-compact operators. We claim that in the case of d-compact
operators the order continuity condition in
Theorem~\ref{t:dPS-main}(\ref{ti:dPS-main:chi-chi}) can be removed.
Indeed, the following result follows immediately from
Proposition~\ref{p:dcomp-set-chi-rho}.

\begin{proposition}
  If $S,T\colon E\to F$ are two operators between Banach lattices such
  that $T$ dominates $S$ and $S$ is d-compact then $\chi(S)\le\chi(T)$.
\end{proposition}

Next, we are going to present several examples of d-topologies.
First, notice that if $E$ has a strong order unit then the d-topology
coincides with the norm topology of $E$.

\begin{example}
  {\it $E=C_0(\Omega)$ where $\Omega$ is a normal topological space.\/}
  The d-convergence in $E$ is exactly the \term{ucc}
  topology, i.e., the topology of uniform convergence on compacta.
  Indeed, let $(x_\alpha)$ be a net in $E$ such that
  $x_\alpha\xrightarrow{\rm ucc}0$, and $y\in E_+$. Fix
  $\varepsilon>0$. Then one can find a compact set $K\subseteq\Omega$
  such that $y(t)\le\varepsilon$ whenever $t\notin K$. There exists an
  index $\alpha_0$ such that $\abs{x_\alpha(t)}\le\varepsilon$
  whenever $t\in K$ and $\alpha\ge\alpha_0$. Then
  $\abs{x_\alpha(t)}\wedge y(t)\le\varepsilon$ for every $t\in\Omega$,
  so that $(x_\alpha)$ is d-null.
  
  Conversely, suppose that $(x_\alpha)$ is d-null in $E$ and $K$ is a
  compact subset of $\Omega$. Let $y\in E$ such that $0\le y\le \one$
  and $y(t)=1$ whenever $t\in K$.  Since $\abs{x_\alpha}\wedge y\to
  0$, it follows that for every $\varepsilon>0$ there exists an index
  $\alpha_0$ such that $\abs{x_\alpha(t)}\wedge y(t)\le\varepsilon$
  whenever $\alpha\ge\alpha_0$. This implies
  $\abs{x_\alpha(t)}\le\varepsilon$ whenever $t\in K$ and
  $\alpha\ge\alpha_0$.
  
  Notice that if $x\in\star{E}$ then $\hat{x}\in E^d$ if and only if
  $x(t)\approx 0$ whenever $t\in\star{K}$ for some compact
  $K\subseteq\Omega$, or, equivalently, if $x(t)\approx 0$ for every
  nearstandard $t\in\star{\Omega}$.
\end{example}

\begin{example}\lab{e:c-0}
  $E=c_0$. We will show that $\widehat{E}\neq\iE\oplus E^d$.
  It follows from the previous example that the
  d-convergence on $c_0$ is exactly coordinate-wise convergence.
  Let $N\in\star{\mathbb N}\setminus{\mathbb N}$ and
  $x=\sum_{k=1}^Ne_k$. Assume that $\hat{x}=\hat{y}+\hat{z}$ for some
  $y,z\in\star{E}$ such that $\hat{y}\in\iE$ and $\hat{z}\in E^d$.
  Then $\norm{\hat{y}-\hat{v}}<\frac{1}{4}$ for some $\hat{v}\in
  I(E)$. Set $D=\bigl\{i\in\star{\mathbb N}\mid
  \abs{z_i}<\frac{1}{2}\bigr\}$. Then $\mathbb N\subseteq D$ because
  $z_i\approx 0$ for each standard $i$. By the Overspill Principle there
  exists $n\in D$ such that $n\notin \mathbb N$ and $n\le N$. Then
  $y_n\approx x_n-z_n\ge\frac{1}{2}$ and, therefore,
  $v_n\ge\frac{1}{4}$. But $v_n$ must be infinitesimal because $v$ is
  dominated by a standard sequence, a contradiction. Thus,
  $\hat{x}\notin\iE\oplus E^d$ so that $\iE\oplus E^d$ is a proper
  subset of $\widehat{E}$. It also follows that the embedding of $E^d$
  into $\check{E}$ given be $\hat{x}\mapsto\check{x}$ is not onto,
  because otherwise we would have $\check{x}=\check{z}$ for some
  $\hat{z}\in E^d$, and this would imply $\hat{x}=\hat{y}+\hat{z}$ for
  $\hat{y}=\hat{x}-\hat{z}\in\iE$.
  
  Furthermore, if $A=\{e_k\}_{k\in\mathbb N}$, it can be easily
  verified that $A$ is d-compact in $c_0$. On the other hand, if we
  consider its convex hull $\co A$ then clearly $x\in\star{\co A}$, so
  that $\co A$ is not relatively d-compact. Thus the convex hull of a
  d-compact set need not be d-compact, and the d-topology need not be
  locally convex.
\end{example}

\begin{remark} \lab{r:E+Ed}
  It follows from Remark~\ref{r:iE-projband} that if $E$ is a
  Banach lattice such that $c_0$ is not lattice finitely representable
  in $E$, then $\widehat{E}=\iE\oplus E^d$. In this case the map
  $\hat{x}\mapsto\check{x}$ is an isometry between $E^d$ and
  $\check{E}$. If, in addition, $E$ is atomic with order continuous norm
  (e.g., $E=\ell_p$, $1\le p<\infty$) then $\iE=E$, so that
  $\widehat{E}=E\oplus E^d$. It follows that every bounded set in $E$
  is relatively d-compact, and every $E$-valued bounded operator is
  d-compact.
\end{remark}

\begin{example} {\it $E=L_p(\mu)$
  where $1\le p<\infty$ and $\mu$ is a finite measure.\/}
  The d-convergence in $E$ is exactly the
  convergence in measure. To show this, let
  $x_\alpha\xrightarrow{\mu}0$ in $E$. Clearly $\abs{x_\alpha}\wedge
  \one$ converges to zero in $L_p$-norm. Similarly,
  $\abs{x_\alpha}\wedge s\to 0$ for every simple function $s\in E$.
  Let $y\in E_+$. For each positive $\varepsilon$ there exists
  a simple function $s\in E$ such that $\norm{s-y}\le\varepsilon$.
  Then:
  \begin{multline*}
    \bignorm{\abs{x_\alpha}\wedge y}\le\bignorm{\abs{x_\alpha}\wedge s}+
    \bignorm{\abs{x_\alpha}\wedge s-\abs{x_\alpha}\wedge y}\\\le
    \bignorm{\abs{x_\alpha}\wedge s}+\norm{s-y}\le
    \bignorm{\abs{x_\alpha}\wedge s}+\varepsilon.
  \end{multline*}
  This yields $\abs{x_\alpha}\wedge y\to 0$, so that $x_\alpha$
  d-converges to zero.
  
  Conversely, suppose $(x_\alpha)$ d-converges to zero in $E$. Fix
  $\varepsilon>0$. Then $\abs{x_\alpha}\wedge\varepsilon\one\to 0$, so
  that we can find $\alpha_0$ such that
  $\bignorm{\abs{x_\alpha}\wedge\varepsilon{\bf
      1}}\le\varepsilon^\frac{p+1}{p}$ whenever $\alpha\ge\alpha_0$.
  This yields:
  \begin{multline*}
    \mu\Bigl(\bigl\{\abs{x_\alpha}>\varepsilon\bigr\}\Bigr)=
    \frac{1}{\varepsilon^p}\int\limits_{\abs{x_\alpha}>\varepsilon}
      \varepsilon^p\,d\mu=
    \frac{1}{\varepsilon^p}\int\limits_{\abs{x_\alpha}>\varepsilon}
      \bigl(\abs{x_\alpha}\wedge\varepsilon)^p\,d\mu\\\le
    \frac{1}{\varepsilon^p}\int
      \bigl(\abs{x_\alpha}\wedge\varepsilon\bigr)^p\,d\mu\le
    \varepsilon
  \end{multline*}
  whenever $\alpha\ge\alpha_0$.
  
  For $x\in\star{E}$ the following lemma guarantees that $\hat{x}\in
  E^d$ if and only if $\abs{x}\wedge\one\approx 0$.
\end{example}

\begin{lemma} \lab{l:wu-Ed}
  If $e$ is a quasi-interior point in a Banach lattice $E$,
  then for $x\in\star{E}$ we have $\hat{x}\in E^d$ if and only if
  $\abs{x}\wedge e\approx 0$.
\end{lemma}

\begin{proof}
  Without loss of generality we assume $x\ge 0$.  It follows from
  Lemma~\ref{l:mud} that $\hat{x}\in E^d$ implies $x\wedge e\approx
  0$. Conversely, assume that $x$ is nearly disjoint with $e$ and let
  $y\in E_+$. By~\cite[Theorem 15.13]{Aliprantis:85} we have
  $\norm{y-y\wedge ne}\to 0$. Fix $\varepsilon>0$. Then there exists
  $n\in\mathbb N$ such that $\norm{y-y\wedge ne}\le\varepsilon$.
  Therefore:
  \begin{multline*}
    \norm{x\wedge y}\le
    \norm{x\wedge y-x\wedge y\wedge ne}+\norm{x\wedge y\wedge ne}\\\le
    \norm{y-y\wedge ne}+\norm{x\wedge ne}\le
    \varepsilon+n\norm{x\wedge e}\le 2\varepsilon.
  \end{multline*}
  It follows that $x$ is nearly disjoint with $y$, so that by
  Lemma~\ref{l:mud} we conclude that $\hat{x}\in{}E^d$.
\end{proof}

The following interesting observation was communicated to the author by
W.B.~John\-son: {\it if $0<q<p<\infty$ then a bounded sequence in
  $L_p(\mu)$ converges to zero in measure
  if and only if it converges to zero in $\norm{\cdot}_q$}. Indeed, it
is well known that convergence in $\norm{\cdot}_q$ implies convergence
in measure. On the other hand, suppose that $(x_n)$ is a norm-bounded
d-null sequence in $L_p(\mu)$ and $0<q<p$, then:
 \begin{multline*}
   \int_{\abs{x_n}<1}\abs{x_n}^q\le
   \biggl(\int_{\abs{x_n}<1}\abs{x_n}^p\biggr)^\frac{q}{p}
   \biggl(\int_{\abs{x_n}<1}{\one}\biggr)^\frac{p-q}{p}\\
   \le\biggl(\int_{\abs{x_n}<1}
     \abs{x_n}^p\wedge{\one}\biggr)^\frac{q}{p}
   \norm{\mu}^\frac{p-q}{p}
   =\bignorm{\abs{x_n}\wedge{\one}}_p^q\cdot
   \norm{\mu}^\frac{p-q}{p}
   \to 0\\
 \end{multline*}
 and
 \begin{displaymath}
   \int_{\abs{x_n}\ge 1}\abs{x_n}^q\le
   \biggl(\int_{\abs{x_n}\ge 1}\abs{x_n}^p\biggr)^\frac{q}{p}
   \biggl(\int_{\abs{x_n}\ge 1}{\one}\biggr)^\frac{p-q}{p}\le
   \norm{x_n}_p^q\cdot\mu\bigl(\abs{x_n}\ge 1\bigr)^\frac{p-q}{p}
   \to 0,
 \end{displaymath}
 so that $\norm{x_n}_q\to 0$. It follows that a bounded subset $A$ in
 $L_p(\mu)$ is (relatively) d-compact if and only if it
 is (relatively) norm compact in $L_q(\mu)$.

\begin{example}
  {\it A positive d-compact operator which is not compact.\/} Consider
  the sequence of intervals
  $A_n=\bigl(\frac{1}{2^n},\frac{1}{2^{n-1}}\bigr]$.  Let
  $f_n=2^{\frac{n}{2}}\chi_{A_n}$. Then $\norm{f_n}_2=1$ for each $n$.
  It is easy to see that the operator $T\colon\ell_2\to L_2[0,1]$
  given by $Te_n=f_n$ is an isometric embedding, hence not compact.
  However, $T$ is d-compact because it is compact as an operator from
  $\ell_2$ to $L_1[0,1]$. Indeed, $T=\sum_{n=1}^\infty e'_n\otimes
  f_n$, but $\sum_{n=1}^\infty\norm{e'_n}_{\ell'_2}\norm{f_n}_1=
  \sum_{n=1}^\infty2^{-\frac{n}{2}}<+\infty$, so that $T$ considered
  as an operator from $\ell_2$ to $L_1[0,1]$ is nuclear, and hence
  compact.
\end{example}

We have mentioned that every order bounded d-compact operator is
AM-compact. The following example was pointed out to the author by
W.B.~Johnson.

\begin{example}
  {\it A d-compact operator which is not AM-compact.\/} Define
  $T\colon L_2[0,1]\to\ell_2$ via $Tx=\bigl(\int
  xr_n\bigr)_{n=0}^\infty$, where $r_n$ is the n-th Rademacher
  function. It follows from Remark~\ref{r:E+Ed} that $T$ is d-compact.
  On the other hand, $Tr_n=e_n$, so that $T$ is not AM-compact.
\end{example}

Recall that a set $A$ in a Banach lattice $E$ is said to be
\term{PL-compact} if it is relatively compact with respect to the
seminorm $f\bigl(\abs{\cdot}\bigr)$ for every $f\in E_+'$.  It was
shown in \cite[Proposition 2.1]{dePagter:88} that $\rho(A)=\chi(A)$
for every PL-compact set $A$ in a Banach lattice with order continuous
norm. Furthermore, if $E'$ and $F$ have order continuous norm then
\cite[Theorem 125.3]{Zaanen:83} guarantees that an order bounded
operator $T\colon E\to F$ is AM-compact if and only if $TB_E$ is
PL-compact. In particular, if $T$ is d-compact then $TB_E$ is
PL-compact.

\begin{example}
  {\it A set which is d-compact but not PL-compact.\/} Consider the
  set $A=\{e_k\}_{k=1}^\infty$ in $\ell_1$. It can be easily seen that
  $A$ is d-compact. Nevertheless $A$ is not PL-compact because if we
  take $f(x)=\sum_{i=1}^\infty x_i$ then $f(\abs{\cdot})$ coincides
  with the norm on $E$, while $A$ is not relatively norm compact.
\end{example}

Finally we would like to mention that it is crucial to describe the
sets that satisfy $\chi(A)=\rho(A)$, because if
$\chi(TB_E)=\rho(TB_E)$ then $\chi(T)=\rho(T)$ and, therefore we can
answer the questions stated in the beginning of the paper in the
affirmative. Indeed, it follows from (\ref{eq:nussbaum}) and
Corollary~\ref{c:rho-of-operator} that
$\re(S)=r(\check{S})\in\sigma(\check{S})\subseteq\sigma(\widetilde{S})=\se(S)$
because $\check{E}$ is a Banach lattice.  Furthermore, if $T$
dominates $S$ then
$\chi(S)=\rho(S)=\norm{\check{S}}\le\norm{\check{T}}=\rho(T)\le\chi(T)$.
Along with (\ref{eq:nussbaum}) this yields $\re(S)\le\re(T)$.

We have $\rho(A)\le\chi(A)$ for every bounded set $A$.  It was already
mentioned that $\rho(A)=\chi(A)$ for every PL-compact set $A$ in a
Banach lattice with order continuous norm and for every d-compact set
in any Banach lattice.  Recall that
\begin{displaymath}
\rho(A)=\max\limits_{\hat{y}\in\widehat{A}}\norm{\check{y}}
\quad\mbox{ and }\quad
\chi(A)=\max\limits_{\hat{y}\in\widehat{A}}\norm{\tilde{y}},
\end{displaymath}
and denote by
$A_\chi$ the set of all the points of
$\widehat{A}$ where the latter maximum is attained.
Clearly, if $\norm{\check{x}}=\norm{\tilde{x}}$ for some $\hat{x}\in
A_\chi$ then $\chi(A)=\rho(A)$.
In particular, $A$ is d-compact then the entire
$\widehat{A}$ is contained in
$E\oplus E^d$ and $\norm{\check{x}}=\norm{\tilde{x}}$ for every
$\hat{x}\in E\oplus E^d$. But clearly d-compactness is a way too
strong condition. It would suffice for just $A_\chi$ and $E\oplus E^d$
to have nonempty intersection. Expressed in standard terms, this idea
gives rise to the following proposition.

\begin{proposition}
  Suppose that $A$ is a bounded set in a Banach lattice $E$. If there
  exists a d-convergent sequence $(x_n)_{n=1}^\infty$ in $A$
  such that $\chi(A)=\chi\bigl(\{x_n\}_{n=1}^\infty\bigr)$
  then $\chi(A)=\rho(A)$.
\end{proposition}

\begin{proof}
  Suppose that $(x_n)_{n=1}^\infty$ d-converges to some
  $x\in E$. Since
  $\chi(A)=\max\limits_{n}\norm{\tilde{x}_n}$
  there exists $n_0\in\star{\mathbb N}$ such that
  $\chi(A)=\norm{\tilde{x}_{n_0}}$. If $n_0\in\mathbb N$ then
  $x_{n_0}\in E$. If $n_0\in\star{\mathbb N}\setminus\mathbb N$
  then Lemma~\ref{l:mud} implies $\hat{x}_{n_0}-\hat{x}\in
  E^d$. In either case $\hat{x}_{n_0}\in E\oplus E^d$ so that
  $\norm{\tilde{x}_{n_0}}=\norm{\check{x}_{n_0}}$, and,
  therefore, $\chi(A)\le\rho(A)$.
\end{proof}

%\bibliographystyle{alpha}
%\bibliography{tv}

\section*{Acknowledgements}

I would like to thank Y.A.~Abramovich, E.Yu.~Emel'\-ya\-nov, C.W.~Hen\-son,
W.B.~John\-son, P.A.~Loeb, and M.P.~Wolff for their interest in the work and
helpful discussions.

\end{document}